\documentclass{article}[11pt]

\usepackage{amssymb,amsthm}

\usepackage{subfigure,amsfonts,amsmath}
\usepackage{epsfig,color}

\newtheorem{theorem}{Theorem}[section]
\newtheorem{lemma}[theorem]{Lemma}
\newtheorem{e-proposition}[theorem]{Proposition}
\newtheorem{corollary}[theorem]{Corollary}
\newtheorem{e-definition}[theorem]{Definition\rm}

\newcommand{\calK}{{\mathcal{K}}}

\newcommand{\jump}[1]{[#1]}

\def\og{\leavevmode\raise.3ex\hbox{$\scriptscriptstyle\langle\!\langle$~}}
\def\fg{\leavevmode\raise.3ex\hbox{~$\!\scriptscriptstyle\,\rangle\!\rangle$}}

\usepackage{epsfig}
\usepackage{psfrag}
\usepackage{amsmath}

\newcommand{\prf}{{\noindent{Proof. }}}
\newcommand{\prfend}{\qed}

\begin{document}

\title{Error estimates for stabilized finite element methods applied to ill-posed problems}

\author{Erik Burman \\Department of Mathematics \\University College
  London}




\maketitle

\begin{abstract}
We propose an analysis for the stabilized finite element methods
proposed in {\emph{E. Burman, Stabilized finite element methods for nonsymmetric, noncoercive, and
  ill-posed problems. Part I: Elliptic equations. SIAM
  J. Sci. Comput., 35(6), 2013}}, valid in the case of ill-posed
problems for which only weak continuous dependence can be assumed. A
priori and a posteriori error estimates are obtained without assuming
coercivity or inf-sup stability of the continuous problem.
\end{abstract}

\section{Introduction}
We are interested in the numerical approximation of ill-posed
problems. Consider as an example the following linear elliptic Cauchy problem.
Let $\Omega$ be a convex polygonal (polyhedral) domain in
$\mathbb{R}^d$ and consider the equation
\begin{equation}
\label{eq:Pb}
\left\{
\begin{array}{rcl}
-\Delta  u & = &  f, \mbox{ in } \Omega\\
u  =  0 \mbox{ and } \nabla u \cdot n &=&\psi \mbox{ on } \Gamma
\end{array}
\right. 
\end{equation}
where $\Gamma \subset \partial \Omega$ denotes a simply connected part of the
boundary and $f \in L^2(\Omega)$, $\psi \in
H^{\frac12}(\Gamma)$. Introducing the spaces
$
V:= \{v \in H^1(\Omega): v\vert_{\Gamma} = 0 \}$ and $W:= \{v \in H^1(\Omega): v\vert_{\Gamma'} = 0 \}
$, where $\Gamma' := \partial \Omega \setminus
\Gamma$ and the forms $a(u,w) = \int_\Omega \nabla u\cdot \nabla w ~\mbox{d}x,$ and $l(w) := 
\int_\Omega f w ~\mbox{d}x + \int_{\Gamma} \psi w ~\mbox{d}s
$ equation \eqref{eq:Pb} may be cast in the abstract weak formulation, find $u \in V$ such that
\begin{equation}\label{abstract_prob}
a(u,w) = l(w)\quad \forall w \in W.
\end{equation} 
It is well known that the Cauchy problem \eqref{eq:Pb} is not
well-posed in the sense of Hadamard. If $l(w)$ is such that a
sufficiently smooth, exact solution exists,
conditional continuous dependence estimates can nevertheless be obtained 
\cite{ARRV09}.

The objective of the present paper is to study numerical methods for ill-posed
problems on the form \eqref{abstract_prob} where 
$a:V \times W \mapsto \mathbb{R}$ and $l:W \mapsto \mathbb{R}$ are a
bilinear and a linear form. 
Assume that the linear form $l(w)$ is such that
the problem
\eqref{abstract_prob} admits a unique solution $u \in V$. Define the following dual norm on $l$,
$
\|l\|_{W'} := \sup_{\substack{w \in W \\\|w\|_{W} = 1}} |l(w)|.
$
{\emph{Observe that we do not assume that \eqref{abstract_prob} admits a
unique solution for all $l(w)$ such that $\|l\|_{W'} < \infty$.}} The stability property we assume to be satisfied by \eqref{abstract_prob} is the following continuous
dependence. \\
\flushleft
{\bf Assumption: continuous dependence on data.} Consider the functional $j:V \mapsto
\mathbb{R}$. Let $\Xi:\mathbb{R}^+\mapsto \mathbb{R}^+$ be a continuous, monotone increasing
function with $\lim_{x \rightarrow 0^+}\Xi(x) = 0$. Assume that for a sufficiently small $\epsilon>0$, there holds
\begin{equation}\label{cont_dep_assump}
\|l\|_{W'} \leq \epsilon \mbox{ in \eqref{abstract_prob} then}
|j(u)| \leq \Xi(\epsilon).
\end{equation}
For the example of the Cauchy problem \eqref{eq:Pb}, it is known \cite[Theorems 1.7
and 1.9]{ARRV09} that if \eqref{eq:Pb} admits a unique solution $u \in
H^1(\Omega)$, a continuous dependence of the form
\eqref{cont_dep_assump}, with $0<\epsilon <1$, holds for 
\begin{multline}\label{stab_1}
\mbox{$j(u):=
\|u\|_{L^2(\omega)}$, $\omega \subset \Omega:\,
\mbox{dist}(\omega, \partial \Omega) =: d_{\omega,\partial \Omega}>0$
with}\\ 
\Xi(x) = C_{u\varsigma} x^\varsigma, \,C_{u\varsigma}>0,\, \varsigma:= \varsigma(d_{\omega,\partial \Omega}) \in (0,1)
\end{multline}
 and for 
\begin{equation}\label{stab_2}
\mbox{$j(u):= \|u\|_{L^2(\Omega)}$ with $\Xi(x) = C_u
(|\log(x)| + C)^{-\varsigma}$ with $C_u, C>0$, $\varsigma \in (0,1)$.}
\end{equation}
Note that to derive these results $l(\cdot)$ is first associated
with its Riesz representant in $W$ (c.f. \cite[equation
(1.31)]{ARRV09} and discussion.) The constant $C_{u\varsigma}
$ in \eqref{stab_1} grows monotonically in $\|u\|_{L^2(\Omega)}$ and
$C_{u}$ in
\eqref{stab_2} grows monotonically in $\|u\|_{H^1(\Omega)}$.
\section{Finite element discretization}
Let $\calK_h$ be a shape
regular, conforming, subdivision of $\Omega$ into non-overlapping 
triangles $\kappa$. The family of meshes $\{\calK_h \}_h$ is
indexed by the mesh parameter $h:=\max(\mbox{diam}(\kappa))$.
Let $\mathcal{F}_I$ be the set of interior faces in $\mathcal{K}_h$
and $\mathcal{F}_\Gamma,\mathcal{F}_{\Gamma'}$ the set of element
faces of $\mathcal{K}_h$ whose interior intersects $\Gamma$ and $\Gamma'$
respectively. We assume that the mesh matches the boundary of $\Gamma$
so that $\mathcal{F}_{\Gamma} \cap \mathcal{F}_{\Gamma'} =
\emptyset$. 
Let $X_h^1$ denote the standard finite element space of continuous, affine
functions. Define
$V_h := V \cap
X_h^1$ and $W_h := W \cap X_h^1$. 
We may then write the finite element method: find $(u_h,z_h) \in V_h \times W_h$ such that,
\begin{equation}\label{stabFEM}
\begin{array}{ll}
\left.\begin{array}{rcl}
a(u_h,w_h) - s_W(z_h,w_h) &=&l(w_h)\\
a(v_h,z_h) + s_V(u_h,v_h) &=& s_V(u,v_h) 
\end{array} \right\}&
\mbox{ for all $(v_h,w_h) \in V_h \times W_h$}.
\end{array}
\end{equation}
A possible choice of stabilization operators for the problem
\eqref{eq:Pb} are 
\begin{equation}\label{stab_choice}
s_V(u_h,v_h):= \sum_{F \in \mathcal{F}_I \cup \mathcal{F}_\Gamma} \int_F h_F
\jump{\partial_n u_h
 } \jump{\partial_n v_h
 }~\mbox{d}s,\quad \mbox{ with } h_F := \mbox{diam}(F)
\end{equation}
and
\begin{equation}\label{stab_adjoint}
 s_W(z_h,w_h):=a(z_h,w_h) \, \mbox{ or } \, s_W(z_h,w_h):= \sum_{F \in \mathcal{F}_I \cup \mathcal{F}_{\Gamma'}} \int_F h_F
\jump{\partial_n z_h
 } \jump{\partial_n w_h
 }~\mbox{d}s
\end{equation}
where $\jump{\partial_n u_h}$ denotes the jump of
 $\nabla u_h \cdot n_F$ for $F \in \mathcal{F}_I$ and when $F
 \in \mathcal{F}_\Gamma$ define $\jump{\partial_n u_h}\vert_F:=\nabla
 u_h \cdot n_{\partial \Omega}$.
Unique existence of $(u_h,z_h)$ solution to \eqref{stabFEM}-\eqref{stab_adjoint} follows using the
arguments of
\cite[Proposition 3.3]{Bu13}. By inspection we have that the system \eqref{stabFEM}
is consistent with \eqref{abstract_prob} for $z_h=0$. Taking the
difference of \eqref{stabFEM} and the relation
\eqref{abstract_prob}, with $w=w_h$, we obtain the Galerkin
orthogonality,
\begin{equation}\label{galortho}
a(u_h - u,w_h) - s_W(z_h,w_h) +
a(v_h,z_h) + s_V(u_h - u,v_h) = 0
\end{equation} 
 for all $(v_h,w_h) \in V_h \times W_h$.
\section{Hypotheses on forms and interpolants}
Consider the general, positive semi-definite,
symmetric stabilization operators,
$
s_V: V_h \times  V_h \mapsto \mathbb{R},\quad s_W: W_h \times W_h \mapsto \mathbb{R}.
$
We assume that $s_V(u,v_h)$, with $u$ the solution of
\eqref{abstract_prob} is explicitly known, it may depend on data from
$l(w)$ or measurements of $u$.
Assume that both $s_V$ and $s_W$ define semi-norms on $H^s(\Omega) +V_h$ and
$H^s(\Omega) +W_h$ respectively, for
some $s\ge1$,
\begin{equation}\label{semi_S}
|v+v_h|_{s_Z} := s_Z(v+v_h,v+v_h)^{\frac12}, \forall v\in
H^s(\Omega),\, v_h \in Z_h, \mbox{ with } Z = V,W.
\end{equation}
Then assume that there exists interpolation operators $i_V: V
\mapsto V_h$ and $i_W: W \mapsto W_h$ and norms $\|\cdot\|_{*,V}$ and $\|\cdot\|_{*,W}$ defined
on $V$ and $W$ respectively, such that the form $a(u,v)$ satisfies the continuities
\begin{equation}
a(v-i_V v, w_h) \leq \|v - i_V v\|_{*,V} |w_h|_{s_W},\, \forall v \in
V,\, w_h \in W_h \label{cont1}
\end{equation}
and for $u$ solution of \eqref{abstract_prob},
\begin{equation}
a(u - u_h, w - i_W w) \leq \delta_l(h) \|w\|_W +
\|w - i_W w\|_{*,W} | u - u_h|_{s_V},\, \forall w \in W. \label{cont2}
\end{equation}
In practice $\delta_l(h)$ only depends on the properties of the
interpolant $i_W$ and the data of the problem (and satisfies $\lim_{h
  \rightarrow 0} \delta_l(h)=0$ provided the data are unperturbed).  We also assume that the interpolants
have the following approximation and stability properties. For sufficiently smooth $v\in V$
there holds, for $t>0$
\begin{equation}\label{approx}
| v-i_V v|_{s_V} + \|v - i_V
v\|_{*,V}\leq C_V(v) h^{t}.
\end{equation}
The factor $C_V(v)>0$ will typically depend on some Sobolev norm of
$v$. For $i_W$ we assume that for some $C_W>0$ there holds
\begin{equation}\label{Wstab}
\|w - i_W
w\|_{*,W}+  |i_W w|_{s_W} \leq   C_W \|w\|_W, \quad \forall w \in W.
\end{equation}

\subsection{Satisfaction of hypothesis for the formulation \eqref{stabFEM} -- \eqref{stab_adjoint}}\label{sec:satisfaction}
Let $i_V$ and $i_W$ be
defined by Scott-Zhang interpolation operators preserving the
Dirichlet boundary conditions. The consistency of $s_V(\cdot,\cdot)$
holds for solutions $u \in H^2(\Omega)$. Consider first the form of
$s_W(\cdot,\cdot)$ in the left definition of \eqref{stab_adjoint}. Define $\|v\|_{*,V} := \|\nabla
v\|_{L^2(\Omega)}$ and $\|w\|_{*,W} := \|h^{-1} w\|_{L^2(\Omega)} + \left(\sum_{F \in \mathcal{F}_I \cup
  \mathcal{F}_\Gamma} h^{-1} \|
w\|^2_{L^2(F)}\right)^{
1/2}$. Using local trace inequalities and the
stability and approximation properties of the Scott-Zhang interpolant we deduce that the
inequalities \eqref{approx}-\eqref{Wstab} hold with $t=1$ and
$C_V(v):=C\|v\|_{H^2(\Omega)}$. The inequality \eqref{cont1} follows by the
Cauchy-Schwarz inequality. To prove \eqref{cont2},
with $\delta(h)=C_W h \|f\|_{L^2(\Omega)}$, integrate by parts in $a(u -
u_h, w- i_W w)$, and use the equation \eqref{eq:Pb}, 
to obtain
\[
a(u -u_h, w- i_W w) = (f, w- i_W w)_{L^2(\Omega)} + \sum_{F \in F_I
  \cup F_{\Gamma}} (\jump{\partial_n (u - u_h)},   w- i_W w)_{L^2(F)}.
\] 
The bound \eqref{cont2} then follows by the Cauchy-Schwarz
inequality, the definitions of $s_V(\cdot,\cdot)$ and
$\|\cdot\|_{*,W}$ and the approximation \eqref{Wstab}. For the variant where $s_W(w_h,z_h):= \sum_{F \in \mathcal{F}_I \cup \mathcal{F}_{\Gamma'}} \int_F h
\jump{\partial_n z_h
 } \jump{\partial_n w_h
 }~\mbox{d}s$ let $\|w\|_{*,V}:=\|h^{-1} w\|_{L^2(\Omega)} +\left(\sum_{F \in \mathcal{F}_I \cup
  \mathcal{F}_{\Gamma'}} h^{-1} \|w\|^2_{L^2(F)}\right)^{
1/2}$ and prove
inequality \eqref{cont1} similarly as \eqref{cont2}
above, but integrating by parts the other way. This latter method has enhanced
 adjoint consistency.

\section{Error analysis}
We will now prove an error analysis using only the continuous
dependence \eqref{cont_dep_assump}. First we prove that assuming 
smoothness of the exact solution the error converges with the
rate $h^t$  in the stabilization semi-norms defined in equation
\eqref{semi_S}. Then we show that the computational error satisfies a perturbation
equation in the form \eqref{abstract_prob}, and that the right hand side of the perturbation equation can be
upper bounded by the stabilization semi-norm. Our error bounds are
then a consequence of the assumption
\eqref{cont_dep_assump}. 
\begin{lemma}\label{lem:stab_conv}
Let $u$ be the solution of
\eqref{abstract_prob} and $(u_h,z_h)$ the solution of
the formulation \eqref{stabFEM} for which
\eqref{semi_S},  \eqref{cont1} and \eqref{approx} hold. Then
$$
|u - u_h|_{s_V} + |z_h|_{s_W} \leq (1+\sqrt{2}) C_V(u) h^{t}.
$$
\end{lemma}
\prf
Let $\xi_h :=
i_V u - u_h$ and write
$
|\xi_h|_{s_V}^2 +  |z_h|^2_{s_W} = s_V(\xi_h,\xi_h) + a(\xi_h,z_h) -
a(\xi_h,z_h) + s_W(z_h,z_h).
$
Using equation \eqref{galortho} we then have
$
|\xi_h|_{s_V}^2 +  |z_h|^2_{s_W} = s_V(i_V u -  u,\xi_h) + a(i_Vu - u,z_h).
$
Applying the Cauchy-Schwarz inequality in the first term of the right
hand side and the continuity \eqref{cont1} in the
second, followed by
\eqref{approx} we may deduce
\[
|\xi_h|_{s_V}^2 +  |z_h|^2_{s_W} \leq |i_V u -  u|_{s_V} |\xi_h|_{s_V}
+  \|i_V u - u\|_{*,V} |z_h|_{s_W} \leq C_V(u) h^{t} (|\xi_h|_{s_V}^2 +  |z_h|^2_{s_W})^{\frac12}.
\]
The claim follows by the triangle
inequality $|u - u_h|_{s_V} \leq |u - i_Vu|_{s_V} +|\xi_h|_{s_V}.$
\prfend
\begin{theorem}\label{thm:cont_dep}
Let $u$ be the solution of
\eqref{abstract_prob} and $(u_h,z_h)$ the solution of the formulation
\eqref{stabFEM} for which \eqref{semi_S}-\eqref{approx} hold. Assume that the problem \eqref{abstract_prob} has the
stability property \eqref{cont_dep_assump}. Then
\begin{equation}\label{aposteriori}
|j(u - u_h)| \leq \Xi(\eta(u_h,z_h))
\end{equation}
where the a posteriori quantity $\eta(u_h,z_h)$ is defined by
$
\eta(u_h,z_h):= \delta_l(h) + C_W(| u - u_h|_{s_V} + |z_h|_{s_W}).
$
For sufficiently smooth $u$ there holds
\begin{equation}\label{apriori}
\eta(u_h,z_h)
 \leq \delta_l(h)  + (1+\sqrt{2})C_W C_V(u) h^{t}.
\end{equation}
\end{theorem}
\prf
Let $e = u - u_h \in V$. By the Galerkin orthogonality there
holds for all $w \in W$
\[
a(e,w) = a(e,w-i_W w) - s_W(z_h,i_W w)=l(w-i_W w) - a(u_h,w - i_W w) - s_W(z_h,i_W w)
\]
and we identify $r \in W'$ such
that $\forall w \in W$,
\begin{equation}\label{pert_right}
(r,w)_{W',W} = l(w-i_W w) - a(u_h,w - i_W w) - s_W(z_h,i_W w).
\end{equation}
We have shown that $e$ satisfies equation \eqref{abstract_prob} with
right hand side $(r,w)_{W',W}$.
Now apply the continuity \eqref{cont2},
Cauchy-Schwarz inequality and the stability \eqref{Wstab} in the right
hand side of \eqref{pert_right} leading to
\[
|(r,w)_{W',W}| = |a(e,w-i_W w) - s_W(z_h,i_W w) |\leq  (\delta_l(h) + C_W |
u - u_h|_{s_V} + C_W  |z_h|_{s_W} )\|w\|_W.
\]
We conclude that 
$
\|r\|_{W'} \leq \delta_l(h) +C_W( | u - u_h|_{s_V} + |z_h|_{s_W}) $
and the claim \eqref{aposteriori} follows by assumption
\eqref{cont_dep_assump}. The upper bound of \eqref{apriori} is a
consequence of Lemma \ref{lem:stab_conv}.
\prfend
\begin{corollary}
Let $u\in H^2(\Omega)$ be the solution of \eqref{eq:Pb} and $u_h,z_h$ the
solution of \eqref{stabFEM}-\eqref{stab_adjoint}. Then the conclusions
of Lemma \ref{lem:stab_conv} and Theorem \ref{thm:cont_dep} hold for
$u-u_h,\,z_h$ with $t=1$ and $j(\cdot)$, $\Xi(\cdot)$
given by \eqref{stab_1} or\eqref{stab_2}. Moreover
$C_{u\varsigma}$ and $C_u$ of \eqref{stab_1} and \eqref{stab_2} are independent of $h$.
\end{corollary}
\prf
In Section \ref{sec:satisfaction} above we showed that the formulation
\eqref{stabFEM}-\eqref{stab_adjoint} satisfies
\eqref{semi_S}-\eqref{approx} and we conclude that Lemma
\ref{lem:stab_conv} and Theorem \ref{thm:cont_dep} hold. For $C_{u\varsigma}$ and $C_u$ of \eqref{stab_1}
and \eqref{stab_2} to be bounded uniformly in $h$, $\|u
- u_h\|_{H^1(\Omega)}$ must be bounded by some constant independent of $h$.
To this end one may prove a discrete Poincar\'e inequality
$\|\nabla u_h\|_{L^2(\Omega)} \leq C_P h^{-1} |u_h|_{s_V}$. Using this
result together with Lemma  \ref{lem:stab_conv} we deduce that $\|\nabla u_h\|_{L^2(\Omega)} \leq
C\|u \|_{H^2(\Omega)}$, which proves the claim.
\prfend
\section{Numerical example}
\begin{figure}
\includegraphics[width=6.5cm]{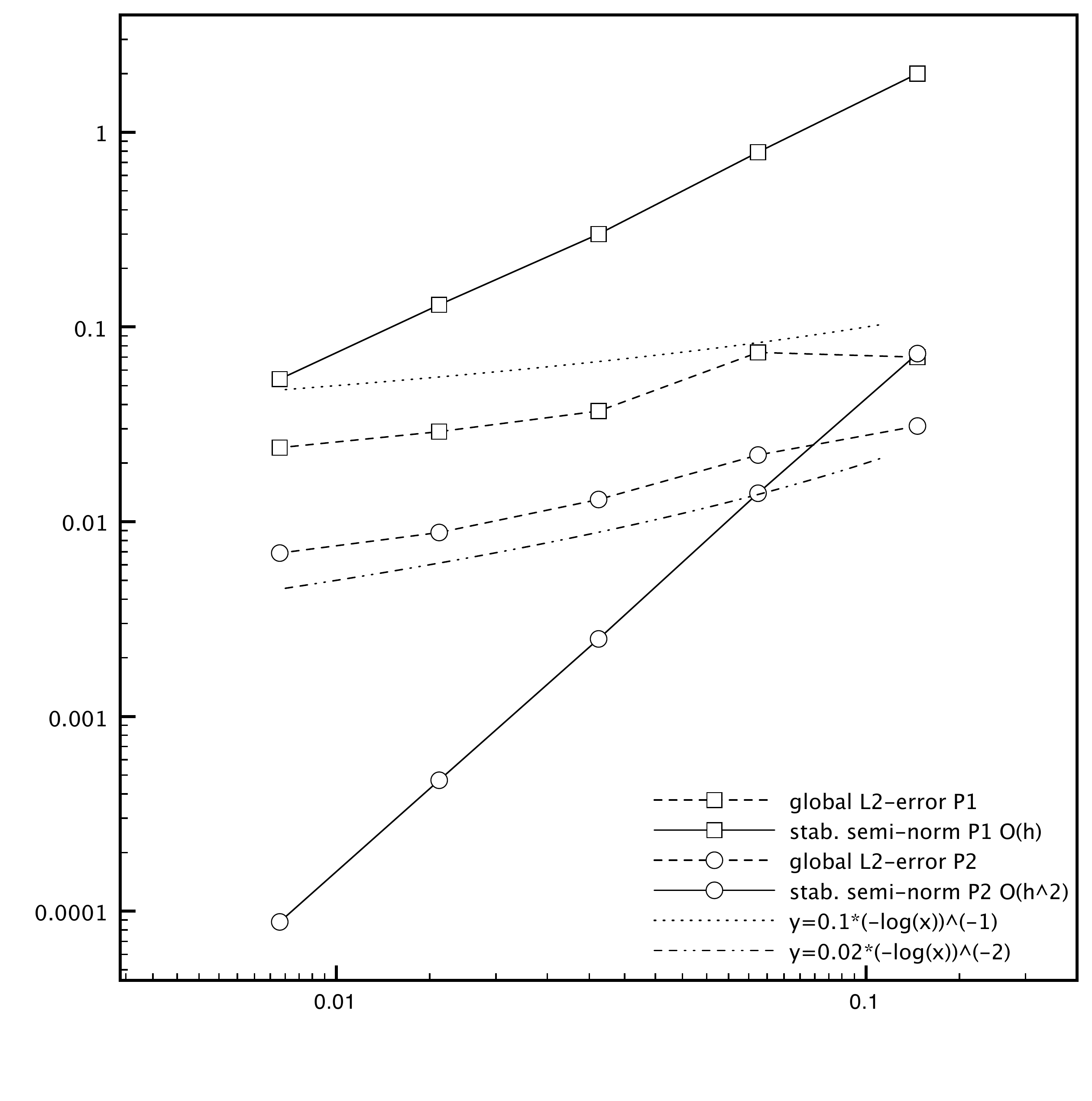}
\includegraphics[width=6.5cm]{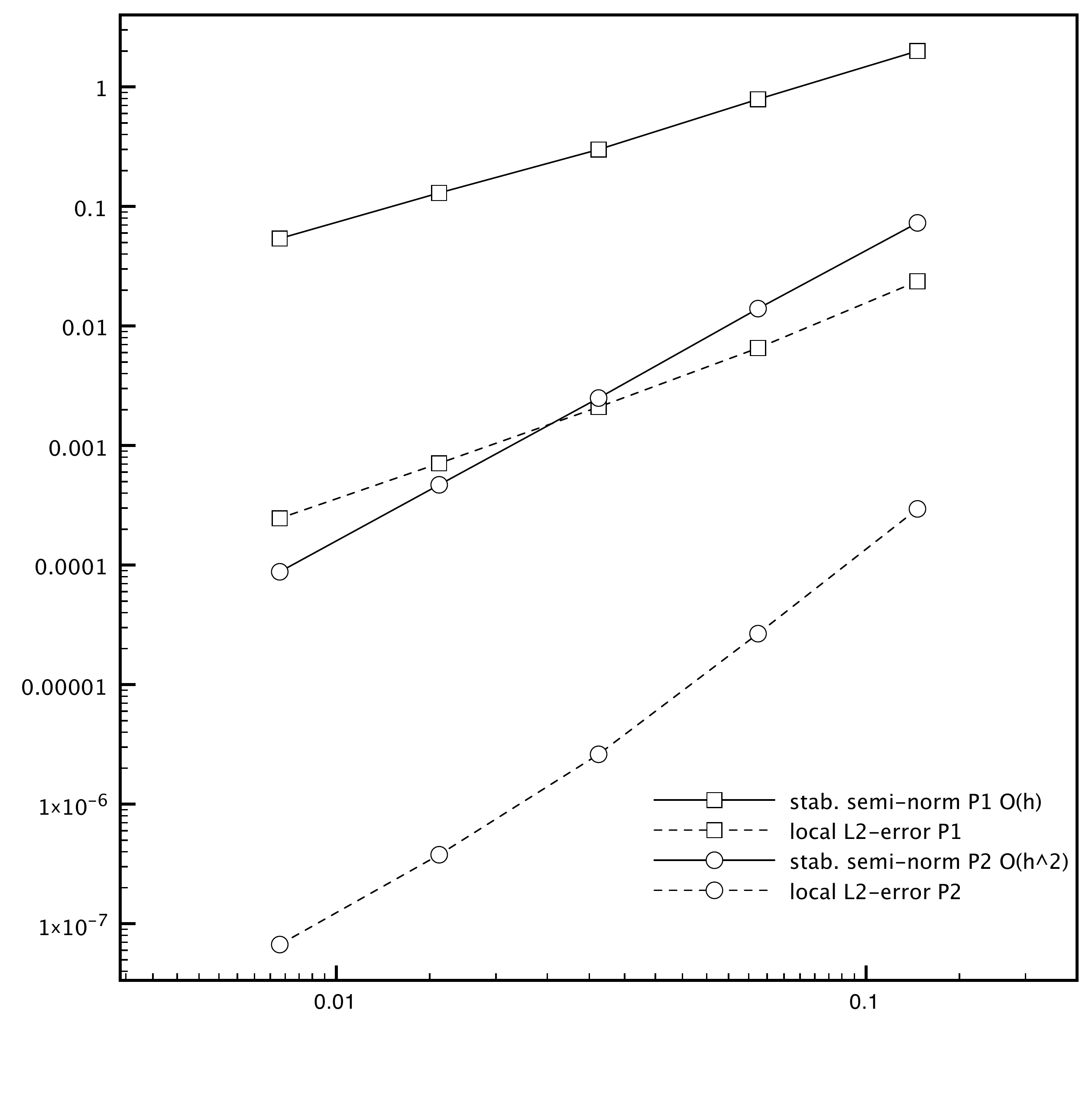}
\caption{Left: convergence of global $L^2$-errors (dashed) and
  stabilization semi-norms (full). Right: convergence
of local $L^2$-errors (dashed) and
  stabilization semi-norms (full). Square markers on curves
  representing $P_1$-approximation and circle markers on curves
  representing $P_2$-approximation.}\label{cauchy_fig}
\end{figure}

\begin{figure}
{\centering{
\includegraphics[width=6.5cm]{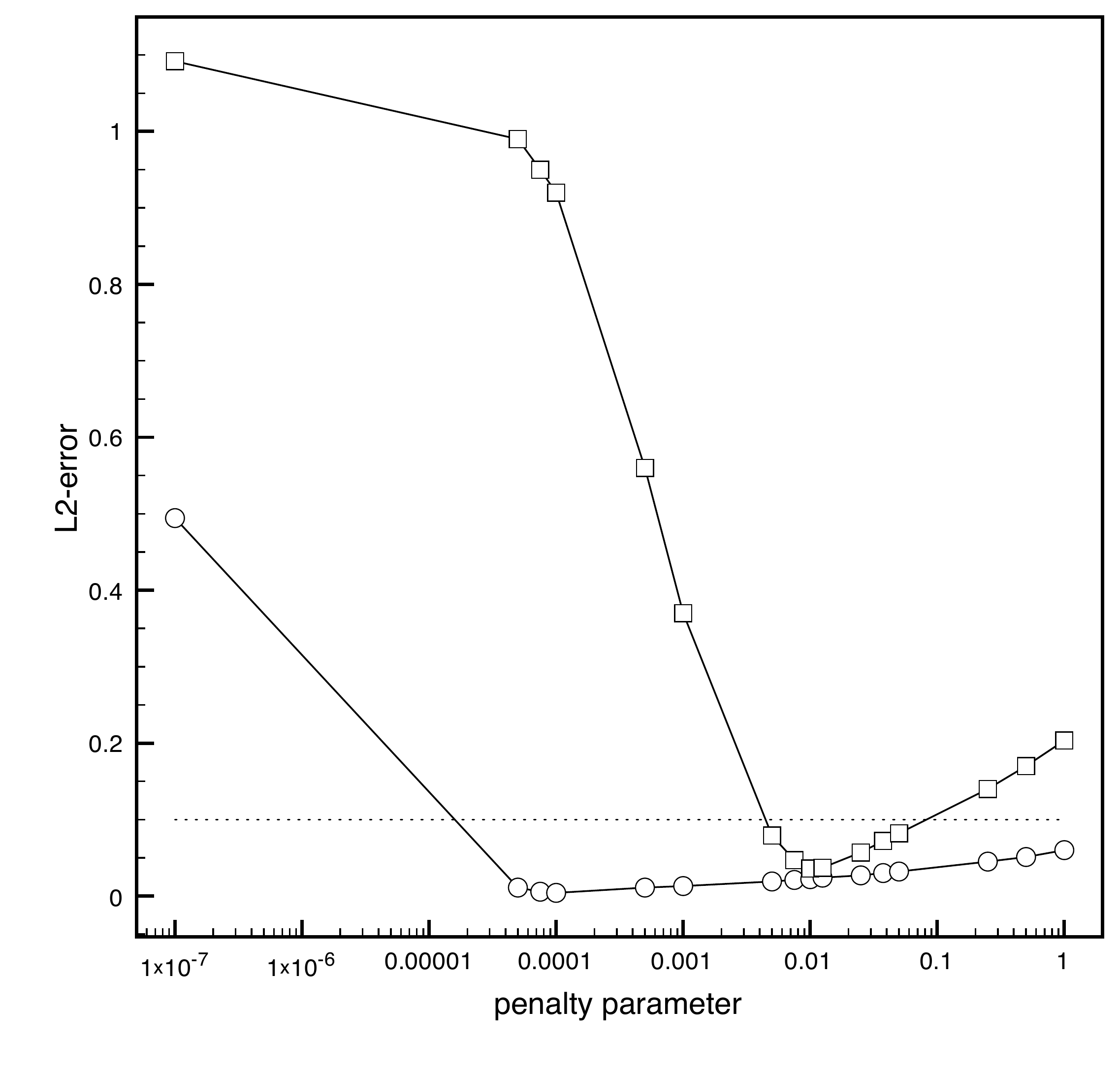}
\caption{
Study of the error under variation of the
parameter $\gamma_V=\gamma_W$. ($P_1$ approximation marked with squares, $P_2$ with circles). }\label{cauchy_param}}}
\end{figure}

To illustrate the theory we recall a numerical example from
\cite{Bu13}. We solve the Cauchy problem \eqref{eq:Pb} on the
unit square $\Omega \in (0,1) \times (0,1)$ with exact solution $u(x,y) = 30x(1-x)y(1-y)$,
$\psi=\nabla u \cdot n_{\partial \Omega}$ and $\Gamma := \{x \in
(0,1), y=0 \}\ \cup \{x=1, \, y \in (0,1) \}$. We compute piecewise affine approximations on a
sequence of unstructured meshes using the method \eqref{stabFEM} and the
stabilizations \eqref{stab_choice} and \eqref{stab_adjoint}$_2$
($\gamma_V = \gamma_W= 0.01$). We also make
a similar series of computations using piecewise quadratic elements
and an added penalty term on the jump of the elementwise Laplacian
following \cite{Bu13} ($\gamma_V = \gamma_W= 0.001$). The results are
reported in Figure \ref{cauchy_fig}. The convergence of the global $L^2$-error and the
stabilization semi-norm is given
in the left plot, compared with theoretically motivated
logarithmic bounds. The local errors in $\omega= (0.5,1) \times (0,0.5)$
are presented in the right plot and we observe that they have $O(h^k)$
convergence where $k$ denotes the polynomial order, similarly as the
stabilization semi-norm. Finally, in Figure \ref{cauchy_param}, we report
a study of the error on a fixed mesh with $64 \times
64$ elements under variation of the penalty parameter in the right plot.
\section{Conclusion and further perspecitives}
Herein we have proposed a framework for the analysis of the 
stabilized methods introduced \cite{Bu13} when applied to ill-posed
problems. The upshot is that error estimates can be obtained using
only continuous dependence properties, without relying on a
well-posedness theory of the continuous problem. 
Important extensions
of the results presented herein are the inclusion of perturbed data and the
exploration of the consequences of adjoint consistency. The latter may allow
for improved estimates, when the error is measured by linear functionals that are in
the range of the adjoint problem. 
\bibliographystyle{abbrv}

\end{document}